\documentclass[a4paper,12pt]{article}
\usepackage{amsmath}
\usepackage{amssymb}
\usepackage{tabularx}
\usepackage{enumerate}
\usepackage{graphicx}
\usepackage{subfigure}
\usepackage{color}
\usepackage[pagewise]{lineno}
\usepackage{longtable}
\usepackage{epstopdf}

\newtheorem{thm}{Theorem}[section]
\newtheorem{lem}[thm]{Lemma}

\newtheorem{exm}{Example}[section]

\newtheorem{rmk}{Remark}[section]
\newtheorem{defi}{Definition}[section]
\newtheorem{pppp}{Proof}

\newcommand{\qed}{\hspace{1em}\mbox{\raisebox{0.65ex}{\fbox{}}}}

\numberwithin{equation}{section}

\newcommand{\be}{\begin{equation}}
\newcommand{\ee}{\end{equation}}
\newcommand\bes{\begin{eqnarray}} \newcommand\ees{\end{eqnarray}}
\newcommand{\bess}{\begin{eqnarray*}}
\newcommand{\eess}{\end{eqnarray*}}

\begin{document}

\thispagestyle{empty}

\title{Effects of impulsive harvesting and an evolving domain in a diffusive logistic model\thanks{The work is partially supported by the NNSF of China (Grant No. 11771381, 61877052).}}

\date{\empty}

\author{Yue Meng$^1$, Zhigui Lin$^1 \thanks{Corresponding author. Email: zglin@yzu.edu.cn (Z. Lin).}$, and Michael Pedersen$^2$\\
{\small 1 School of Mathematical Science, Yangzhou University, Yangzhou 225002, China}\\
{\small 2 Department of Applied Mathematics and Computer Science,}\\
{\small Technical University of Denmark, DK 2800, Lyngby, Denmark}
}

 \maketitle
\begin{quote}
\noindent
{\bf Abstract.} { 
 In order to understand how the combination of domain evolution and impulsive harvesting affect  the dynamics of a population, we propose a diffusive logistic population model with impulsive harvesting on a periodically evolving domain. Initially the ecological reproduction index of the impulsive problem is introduced and given by an explicit formula, which depends on the domain evolution rate and the   impulsive function. Then the threshold dynamics of the population under monotone or nonmonotone impulsive harvesting is established based on this index. Finally numerical simulations are carried out to illustrate our theoretical results, and  reveal that a large domain evolution rate can improve the population survival, no matter which impulsive harvesting takes place. Contrary, impulsive harvesting has a negative effect on the population survival, and can even lead to the extinction of the population.
}


\end{quote}

\section{Introduction}
Mathematical models,  in which reaction diffusion models often focus on the spread and persistence of a population \cite{cc},
have been widely used to study ecological phenomena \cite{mjd}. In nature, species grow and move randomly like e.g. migration of animals,
expansion of invaders and so on. Fish and large mammals give birth in a specific period, so the population experience a birth impulsive growth \cite{ch}.
On the other hand, in order to promote a sustainable development of ecosystems or obtain sustainable economic benefits, population
systems are often disturbed by human management and development strategies \cite{ccw}, such as planting and harvesting, vaccination
or the release of natural enemies. Such stages of disturbance typically occur for a short period of time,
but have significant impacts on the population density and numbers. Classical differential equations are not well suited to describe such phenomena,
where important drivers are non-continuous processes. Hence impulsive differential equations have been introduced to model
and characterize these hybrid discrete-continuous processes.

The theoretical research of impulsive differential equations began by work of   Mil'man and Myshkis in the 1960s \cite{mvd},
and has further developed since  the 80s. An impulsive differential system is a special dynamical system, which has characteristics of both continuous
systems and discrete systems. Its remarkable characteristic is that it can induce short-term rapid changes on the state of systems,
depending on time and state. Impulsive ordinary differential equations have been widely studied in population ecological dynamics,
such as predator-prey systems \cite{ssy}, pest management systems \cite{tsyO}, systems with control strategies \cite{ghj} and so on.
Since drugs are often introduced into the body in the form of pulses by oral or injection in the treatment of diseases, impulsive
ordinary differential equations have also been used to analyze the dynamics of infectious diseases, see \cite{dls, gs, srq} and references therein.

Recently, due to the complex dynamics generated by pulses, it has become an interesting focus of attention how pulses affect a reaction diffusion system.
Some scholars have studied the effects of various responses on the dynamics of impulsive reaction diffusion predator-prey systems \cite{lzj}
and the impacts of seasonal variation on their pattern formation \cite{wxy}. Lewis and Li \cite{lma} investigated a reaction
diffusion model with a seasonal birth pulse and explored how the birth pulse affects the dynamics of the population, including
spreading speed, minimal domain size, travelling waves and complex bifurcations. Considering a nonlocal dispersal stage based on
the work in \cite{lma}, Wu and Zhao \cite{wrw} established threshold-type dynamics of the system when the domain is bounded and
proved the existence of  a spreading speed in an unbounded domain.  Liang, et al. \cite{ljh} extended the growth population model
of pest with multiple pulses, including birth and pesticide applications. They analyzed effects of the birth rate and the killing efficacy
on a traveling wave and spreading speed, and also investigated the optimal timing of use of  pesticides, pest control strategy and spatial
dependent killing rates by numerical simulations.

Previous studies on impulsive differential equations have been carried out on fixed domains. Space is identified as a functional property
of ecological processes in nature \cite{td}, and it does play an important role in the dispersal of the population. For example,
habitat structure can enhance the persistence of the species \cite{esp}. It has been well documented that a spatially inhomogeneous environment affects the population dynamics, disease transmission and species evolution \cite{ly}. With the increased focus  on the role of habitats in ecology,
the influence of shifting domains on species dynamics has become a concern recently, due to the obvious fact that habitats are not constant in nature.
Free boundary problems, modeling situations where shifting boundaries of the habitat are unknown and affected by the movement
of the population, are used to describe invasive species \cite{dyh, lcx} and the spread of infectious diseases \cite{lzg}.
On the other hand, it is well-known that habitats where species live on, often present cyclic variation, which is influenced
by environment factors such as rainfall, temperature and so on. For instance, the increase of rainfall can reduce the riparian
habitat area of waterbirds, in cases where rainfall shows a long-term cyclic behavior, like on the Pampas region of Argentina \cite{cad}.
Such domains, where the changing boundaries are known, are called periodically evolving domains.

Motivated by the above, in order to explore how an evolving domain and impulsive harvesting affect the dynamics of a population,
we consider a diffusive logistic model with impulsive harvesting on a periodically evolving domain. This model describes the situation where a harvesting pulse,
which is described by a function $g$, takes place at every time $nT(n=0,1,2,...)$ during the continuous growth and dispersal process of a population.
The density of the population at the end of the harvesting stage is given by the function $g$ applied to the density of the population at the beginning of the stage.
 Since $0<g(u)/u<1$ we have that $1-g(u)/u$ denotes harvesting rate. During the dispersal stage, the population diffuses by the coefficient $d(>0)$
 and follows a logistic equation. $\alpha(>0)$ represents the   intrinsic growth rate of the population and the effect of the interspecific
 competition among individuals is denoted by $\gamma(>0)$.  Let $u(t,x)$ be the density of the population at time $t$ and point $x\in[0,l]$
 with initial density $u(0,x)=u_0(x)$. With a Dirichlet boundary condition, the population model in a one-dimensional space is introduced as follows:
\begin{eqnarray}
\left\{
\begin{array}{lll}
\frac{\partial u}{\partial t}=d\frac{\partial^2 u}{\partial x^2}+\alpha u-\gamma u^2,\; &\,x\in(0,l),\ t\in((nT)^+,(n+1)T],\ n=0,1,2,..., \\[2mm]
u(t,0)=u(t,l)=0,\; &\,  t>0,\\[2mm]
u(0,x)=u_0(x)\geq,\not\equiv0,\; &\, x\in[0,l],\\[2mm]
u((nT)^+,x)=g(u(nT,x)),\; &\,x\in(0,l),\ n=0,1,2,...
\end{array} \right.
\label{a01}
\end{eqnarray}
In this paper, we make the following assumption about the pulse function $g$:

\medskip
(A1) \ \textit{ $g(u)$ is a once continuously differentiable function for $u\geq0$, $g(0)=0$, $g'(0)>0$, and for $u>0$, $g(u)>0$,
$g(u)/u$ is nonincreasing with respect to $u$ and $0<g(u)/u<1$.}

The pulse functions $g$ satisfying (A1) usually takes the form ( including the Beverton-Holt function):
\bes
g(u)=\frac{mu}{a+u}
\label{a02}
\ees
with $m>0$ and $a>0$ as in \cite{bh}, or the Ricker function:
\bes
g(u)=ue^{r-bu}
\label{a03}
\ees
with $r>0$ and $b>0$ as in \cite{rick}.

Model (\ref{a01}) can also be used to describe that the population grows and diffuses following the logistic equation outside winter,
while during winter, the population stops birth and movement and rely on individuals to survive the winter, for example,
Aedes aegypti mosquitoes live and lay eggs in soils near water during warm seasons, and during winter they stop moving,
die in large numbers and some of them eventually survive in the form of overwintering eggs \cite{csr}.

As in \cite{ej} and \cite{plq}, let $(0,l(t))$ be a periodically evolving domain at time $t>0$ with the shifting boundary
$l(t).$ $u(t,x(t))$ denotes the population density at time $t>0$ and location $x(t)\in (0,l(t))$.

Let $x(t)=m(t)$ and $x(t)=n(t)(m(t)\leq n(t))$ be two arbitrary end points of the interval where $x(t)$ varies from 0 to $l(t)$.
Then, $(m(t),n(t))\subset (0,l(t))$ denotes an evolving domain at time $t>0$ with shifting starting point $m(t)$ and ending point $n(t)$.
By the principle of mass conservation we have that
\bes
\begin{array}{llll}
\frac{d}{dt}\int_{m(t)}^{n(t)}u(t,x(t))dx&=&du_x(t,n(t))-du_x(t,m(t))+\int_{m(t)}^{n(t)}f(t,u(x(t)))dx,\\[2mm]
&=&\int_{m(t)}^{n(t)}(du_{xx}+f(u))dx,
\end{array}
\label{a04}
\ees
where $f(u)=\alpha u-\gamma u^2$. We further assume that the evolution of the domain is uniform and isotropic, in other words,
the domain evolves by the same ratio in all directions as time increases. One possibility can be described as
\bes
\begin{array}{llll}
x(t)=\rho(t)y, \ \ y\geq0,
\end{array}
\label{a05}
\ees
where the positive continuous function $\rho(t)$ is called evolution rate. Here, because of the periodic evolution of the domain $(0,l(t))$,
$\rho(t)$ is T-periodic in time i.e.
\bes
\begin{array}{llll}
\rho(t+T)=\rho(t)
\end{array}
\label{a06}
\ees
for some $T>0$ and $\rho(0)=1$. A study of periodically evolving domains can be found in e.g. \cite{jdh, pr}. If $\dot{\rho}(t)\geq0$,
the domain is called a growing domain, see \cite{ma} and references therein, and a shrinking domain $(\dot{\rho}(t)\leq0)$
has been discussed in reference \cite{wn} and others.

Let $l(0)=l_0$, then the evolving domain can be written $(0,l(t))=(0,\rho(t)l_0)$. Set $m(t)=\rho(t)y_1$ and $n(t)=\rho(t)y_2$,
where $y_1,y_2\in (0,l_0)$. Noting that $u(t, x(t))=u(t, \rho(t)y)$, we define $u(t, x(t))=v(t, y)$, which together with (\ref{a05}) yields
\bes
u_{x}=\frac{1}{\rho(t)}v_{y}, \ \ u_{xx}=\frac{1}{\rho^2(t)}v_{yy}.
\label{a07}
\ees
From the left hand side of (\ref{a04}), we have
\bes
\begin{array}{llll}
\frac{d}{dt}\int_{m(t)}^{n(t)}u(t, x(t))dx&=&\frac{d}{dt}\int_{\rho(t)y_1}^{\rho(t)y_2}u(t, \rho(t)y)d(\rho(t)y),\\[2mm]
&=&\frac{d}{dt}\int_{y_1}^{y_2}v(t, y)\rho(t)dy,\\[2mm]
&=&\int_{y_1}^{y_2}(v_t(t, y)\rho(t)+v(t, y)\dot\rho(t))dy.
\end{array}
\label{a08}
\ees
Using (\ref{a07}), the right hand side becomes
\bes
\begin{array}{llll}
\int_{m(t)}^{n(t)}(du_{xx}+f(u) )dx=\int_{y_1}^{y_2}(\frac{d}{\rho^2(t)}v_{yy}+f(v(t, y)))\rho(t)dy.
\end{array}
\label{a09}
\ees
From (\ref{a08}) and (\ref{a09}),  equation (\ref{a04}) can be written as follows
\bes
\int_{y_1}^{y_2}(v_t+\frac{\dot{\rho}(t)}{\rho(t)}v)dy=\int_{y_1}^{y_2}(\frac{d}{\rho^2(t)}v_{yy}
 + f(v(t, y)))dy,  & \, t>0.
\label{a10}
\ees
Since $m(t)$ and $n(t)$ are arbitrary, then (\ref{a10}) holds for any $y_1,y_2\in(0,l_0)$.
Therefore, (\ref{a01}) is transformed into the following problem on a fixed domain
\begin{eqnarray}
\left\{
\begin{array}{lll}
v_t=\frac{d}{\rho^2(t)}v_{yy}+(\alpha-\frac{\dot{\rho}(t)}{\rho(t)}) v-\gamma v^2,\; &\,y\in(0,l_0),\ t\in((nT)^+,(n+1)T], \\[2mm]
v(t,0)=v(t,l_0)=0,\; &\,  t>0,\\[2mm]
v(0,y)=v_0(y)\geq,\not\equiv0,\; &\, y\in(0,l_0),\\[2mm]
v((nT)^+,y)=g(v(nT,y)),\; &\,y\in(0,l_0),\ n=0,1,2,...
\end{array} \right.
\label{a11}
\end{eqnarray}

When $g(v)=v$, that is, the case when impulsive harvesting does not occur,  (\ref{a11}) is reduced to a classical logistic
equation on a periodically evolving domain, previously discussed by Jiang and Wang \cite{jdh}. They proved persistence
and extinction of species based on the critical value $D_0:=\frac{aT}{\lambda_1 \int_0^T \rho^{-2}(t)dt}$,
where $\lambda_1(>0)$ is the principal eigenvalue of $-\Delta$ in $(0,l_0)$ under the Dirichlet boundary condition.
The results indicate that species become extinct when $d\in(D_0,+\infty)$, on the other hand, species persist if $d\in (0,D_0)$.
They also analyzed effects of the evolution rate on the persistence of species.

We are interested in what amusing changes the impulsive harvesting will impose on the system (\ref{a11}),
whether a new threshold value like $D_0$ in \cite{jdh} in terms of the pulse function can be introduced to establish
threshold-type results for the persistence and extinction of the population, and whether these results are consistent
with results in \cite{jdh}. This paper also aims to explore how regional evolution affects the dynamic behaviour of
the population when impulsive harvesting takes place, and is organized as follows: In the next section,
a new threshold value is introduced, which is the ecological reproduction index of the problem with impulses.
Moreover, we provide an explicit formula of the ecological reproduction index and analyze the relationship between
domain evolution and this index. Section 3 is concerned with threshold-type results on the asymptotic behaviours of the
solution to problem (\ref{a11}) when the impulsive harvesting is either monotone or not. Numerical simulations are
carried out to understand the effects of the evolution rate and impulsive harvesting on the dynamics of the population,
and then biological explanations are given in the Section 4. The last section contains a discussions of the results.

\section{The ecological reproduction index}
In this section, a new threshold value, which is the ecological reproduction index of the problem with impulses, is introduced and analyzed.

Linearizing problem (\ref{a11}) at $v=0$, which represents the disease-free equilibrium epidemiologically,  we obtain the following linearized system
\begin{eqnarray}
\left\{
\begin{array}{lll}
v_t=\frac{d}{\rho^2(t)}v_{yy}+(\alpha-\frac{\dot{\rho}(t)}{\rho(t)})v, &y\in(0,l_0),\ t\in((nT)^+,(n+1)T],\ n=0,1,2,..., \\[2mm]
v(t,0)=v(t,l_0)=0,\; &\,  t>0,\\[2mm]
v((nT)^+,y)=g'(0)v(nT,y),\; &\,y\in(0,l_0),\ n=0,1,2,...
\end{array} \right.
\label{b01}
\end{eqnarray}
We first consider the auxiliary system as follows
\begin{eqnarray}
\left\{
\begin{array}{lll}
v_t=(\alpha-\frac{\dot{\rho}(t)}{\rho(t)})v,\; &\,t\in((nT)^+,(n+1)T],\ n=0,1,2,..., \\[2mm]
v((nT)^+)=g'(0)v(nT),\; &\,n=0,1,2,...
\end{array} \right.
\label{b02}
\end{eqnarray}
As in \cite{bzg}, let $E(t,s)$ be the evolution operator of the problem
\begin{eqnarray}
\left\{
\begin{array}{lll}
v_t=-\frac{\dot{\rho}(t)}{\rho(t)}v,\; &\,t\in((nT)^+,(n+1)T],\ n=0,1,2,..., \\[2mm]
v((nT)^+)=g'(0)v(nT),\; &\,n=0,1,2,...
\end{array} \right.
\label{b03}
\end{eqnarray}
Then, the evolution operator $E(t,s)$ can be written as
\[
E(t,s)=e^{-\int_s^t \frac{\dot{\rho}(\tau)}{\rho(\tau)}d\tau}(g'(0))^k,
\]
where $k$ represents the number of impulsive points on $[s,t)$. Since $\int_s^t \frac{\dot{\rho}(\tau)}{\rho(\tau)}d\tau $ is bounded for any $t>s$,
there exists a positive constant $K$ such that
\[
\parallel E(t,s) \parallel \leq K,\ t\geq s,\ s\in \mathbb{R}.
\]

Let $C_T$ be the Banach space given by
{\footnotesize $$ C_T=\left\{\omega\mid
\omega\in C((nT,(n+1)T]),\, \omega(t+T)=\omega(t)\, \textrm{for}\, t\in \mathbb{R}, \, \omega((nT)^+)=\omega(((n+1)T)^+),\, n\in \mathbb{Z}
 \right\}
$$}
with the  norm $\parallel \omega \parallel=\sup\limits_{t\in[0,T]}|\omega(t)|$ and the positive cone
$C_T^+:=\{ \omega\in C_T\mid \omega(t)\geq 0,\ \forall t\in \mathbb{R} \}$. A linear operator on $C_T $ can be introduced by
\[
[L\omega](t)=\int_0^\infty \alpha E(t,t-s) \omega(t-s)ds.
\]
It is clear that $L$ is positive and compact on $C_T$. According to \cite{bzg}, we define the spectral radius of $L$
\[
\Re_0=r(L)
\]
as the basic reproduction number of the periodic impulsive system (\ref{b02}). Moreover, Theorem 2 in \cite{bzg} shows that $\Re_0=\mu$,
where $\mu$ satisfies the following problem
\[
\left\{
\begin{array}{lll}
v_t=(\frac{\alpha}{\mu}-\frac{\dot{\rho}(t)}{\rho(t)})v,\; &\,t\in((nT)^+,(n+1)T],\ n=0,1,2,..., \\[2mm]
v((nT)^+)=g'(0)v(nT),\; &\,n=0,1,2,...,\\[2mm]
v(0)=v(T).
\end{array} \right.
\]

Now we consider the following periodic eigenvalue problem:
\begin{eqnarray}
\left\{
\begin{array}{lll}
\phi_t=\frac{d}{\rho^2(t)}\phi_{yy}+(\frac{\alpha}{R_0}-\frac{\dot{\rho}(t)}{\rho(t)})\phi,\; &\,y\in(0,l_0),\ t\in((nT)^+,(n+1)T],\ n=0,1,2,..., \\[2mm]
\phi(t,0)=\phi(t,l_0)=0,\; &\,  t>0,\\[2mm]
\phi(0,y)=\phi(T,y),\;&\,y\in[0,l_0],\\[2mm]
\phi((nT)^+,y)=g'(0)\phi(nT,y),\; &\,y\in(0,l_0),\ n=0,1,2,...
\end{array} \right.
\label{b04}
\end{eqnarray}
The theory for the basic reproduction number for the impulsive reaction-diffusion system is not fully established, but similar to the above,
we introduce the ecological reproduction index $R_0$ for our impulsive system by solving problem (\ref{b04}), which can give the explicit
formula of $R_0$ and find a corresponding positive eigenfunction $\phi(t,y)$.

\begin{thm}
The ecological reproduction index of problem $(\ref{a11})$ can be explicitly expressed as
\bes
\begin{array}{llllll}
R_0=\frac{\alpha}{\frac{d\lambda_1}{T}\int_0^T\frac{1}{\rho^2(t)}dt-\frac{1}{T}\ln g'(0)},
\end{array}
\label{b05}
\ees
where $\lambda_1(>0)$ is the principal eigenvalue of $-\Delta$ in $(0,l_0)$ under Dirichlet boundary condition, and $R_0$ is monotonically increasing with respect to $\rho(t)$.
\end{thm}

\textbf{Proof:}
Let
\[\phi(y,t)=f(t)\psi(y),
\]
where $\psi(y)$ is the eigenfunction related to $\lambda_1$ in the eigenvalue problem
\begin{eqnarray}
\left\{
\begin{array}{lll}
-\psi''=\lambda_1\psi, \; &\,y\in(0,l_0), \\[2mm]
\psi(0)=\psi(l_0)=0,\; &\,
\end{array} \right.
\label{b06}
\end{eqnarray}
then problem (\ref{b04}) becomes
\small{\begin{eqnarray}
\left\{
\begin{array}{lll}
f'(t)\psi(y)=\frac{d}{\rho^2(t)}f(t)\psi_{yy}(y)+(\frac{\alpha}{R_0}-\frac{\dot{\rho}(t)}{\rho(t)})f(t)\psi(y), \; &\,y\in(0,l_0),\ t\in((nT)^+,(n+1)T], \\[2mm]
\psi(0)=\psi(l_0)=0,\\[2mm]
f(0)=f(T),\\[2mm]
f((nT)^+)=g'(0)f(nT),\; &\, n=0,1,2,...
\end{array} \right.
\label{b07}
\end{eqnarray}}
By separating variables, it follows from the first equation of problem (\ref{b07}) that
\[
\frac{f'(t)+(\frac{\dot{\rho}(t)}{\rho(t)}-\frac{\alpha}{R_0})f(t)}{\frac{d}{\rho^2(t)}f(t)}=\frac{\psi''(y)}{\psi(y)}=-\lambda_1.
\]
Then, the first equation in problem (\ref{b07}) becomes
\bes
\begin{array}{llll}
f'(t)+(\frac{\dot{\rho}(t)}{\rho(t)}-\frac{\alpha}{R_0}+\frac{d\lambda_1}{\rho^2(t)})f(t)=0.
\end{array}
\label{b08}
\ees
By solving equation (\ref{b08}), we find
\[
f(t)=Ce^{\int_0^t -\frac{\dot{\rho}(\tau)}{\rho(\tau)}+\frac{\alpha}{R_0}-\frac{d\lambda_1}{\rho^2(\tau)}d\tau},\ t\in(0^+,T],
\]
where the initial value $C$ satisfies $C=f(0^+)=g'(0)f(0)$.
Then,
\[
f(T)=g'(0)f(0)e^{\int_0^T -\frac{\dot{\rho}(\tau)}{\rho(\tau)}+\frac{\alpha}{R_0}-\frac{d\lambda_1}{\rho^2(\tau)}d\tau},
\]
which together with (\ref{a06}) and the third equation in problem (\ref{b07}) yields
\[
g'(0)e^{\frac{\alpha T}{R_0}-d\lambda_1\int_0^T\frac{1}{\rho^2(t)}dt}=1.
\]
Therefore,
\[
R_0=\frac{\alpha}{\frac{d\lambda_1}{T}\int_0^T\frac{1}{\rho^2(t)}dt-\frac{1}{T}\ln g'(0)}.
\]
From the explicit formula (\ref{b05}), the monotonicity is easily obtained. \qed\hfill

We note that when $g(v)=v$, then $R_0$ becomes
\bes
R:=\frac{\alpha T}{d\lambda_1\int_0^T\frac{1}{\rho^2(t)}dt},
\label{b09}
\ees
which is consistent with the threshold value in problem without impulses discussed in \cite{jdh}.

\begin{rmk}
When $g'(0)\leq 1$, it is obvious that the ecological reproduction index $R_0>0$. If $g'(0)>1$, the positivity of $R_0$ can not be guaranteed,
and  then $R_0$ as an ecological reproduction index makes no sense. In this situation we can consider the following problem:
\begin{eqnarray}
\left\{
\begin{array}{lll}
v_t=\frac{d}{\rho^2(t)}v_{yy}+(\alpha+M-\frac{\dot{\rho}(t)}{\rho(t)}) v-\gamma v^2-Mv,\; &\,y\in(0,l_0),\ t\in((nT)^+,(n+1)T], \\[2mm]
v(t,0)=v(t,l_0)=0,\; &\,  t>0,\\[2mm]
v(0,y)=v_0(y)\geq,\not\equiv0,\; &\, y\in(0,l_0),\\[2mm]
v((nT)^+,y)=g(v(nT,y)),\; &\,y\in(0,l_0),\ n=0,1,2,...,
\end{array} \right.
\label{b10}
\end{eqnarray}
which is equivalent to problem \eqref{a11}.
Then, the corresponding periodic eigenvalue problem is as follows
\begin{eqnarray}
\left\{
\begin{array}{lll}
\phi_t=\frac{d}{\rho^2(t)}\phi_{yy}+(\frac{\alpha+M}{R_0^*}-\frac{\dot{\rho}(t)}{\rho(t)})\phi-M\phi,\; &\,y\in(0,l_0),\ t\in((nT)^+,(n+1)T], \\[2mm]
\phi(t,0)=\phi(t,l_0)=0,\; &\,  t>0,\\[2mm]
\phi(0,y)=\phi(T,y),\;&\,y\in[0,l_0],\\[2mm]
\phi((nT)^+,y)=g'(0)\phi(nT,y),\; &\,y\in(0,l_0),\ n=0,1,2,...
\end{array} \right.
\label{b11}
\end{eqnarray}
Therefore, the ecological reproduction index $R_0^*:=\frac{\alpha+M}{\frac{d\lambda_1}{T}\int_0^T\frac{1}{\rho^2(t)}dt-\frac{1}{T}\ln g'(0)+M}$,
where $M=\frac{1}{T}|\ln g'(0)|$ can be chosen to make sure that  $R_0^*>0$. The threshold-type dynamics of problem \eqref{a11},
which is established based on the new ecological reproduction index $R_0^*$, is the same as if based on $R_0$.
\end{rmk}

The following result holds for reaction diffusion problems without impulses \cite{lx, zxq}, and also holds for our impulsive problem:
\begin{lem}
${\rm sign}(R_0-1)={\rm sign}\lambda^*$, where $\lambda^*$ satisfies
\[
\left\{
\begin{array}{lll}
\phi_t=\frac{d}{\rho^2(t)}\phi_{yy}+(\alpha-\frac{\dot{\rho}(t)}{\rho(t)})\phi-\lambda^*\phi,\; &\,y\in(0,l_0),\ t\in((nT)^+,(n+1)T],\ n=0,1,2,..., \\[2mm]
\phi(t,0)=\phi(t,l_0)=0,\; &\,  t>0,\\[2mm]
\phi((nT)^+,y)=g'(0)\phi(nT,y),\; &\,y\in(0,l_0),\ n=0,1,2,...
\end{array} \right.
\]
\end{lem}
In fact, $\lambda^*=\frac{1}{T}\ln g'(0)+\alpha-\frac{d\lambda_1}{T}\int_0^T\frac{1}{\rho^2(t)}dt$. The result is obvious.

\section{Asymptotic behaviors of the solution}
In this section, we first explore the asymptotic behaviors of the solution to problem (\ref{a11}) based on the threshold value $R_0$ in
the case where $g$ is monotone, and then we consider the nonmonotone case.
\subsection{Monotone Case}
We start by stating the assumptions:

\medskip
(B1) \textit{$g(u)$ is nondecreasing for $u\geq0$.}

\medskip
(B2) \textit{There are positive constants $D$, $\sigma<\theta$, and $\nu>1$ such that $g(u)\geq g'(0)u-Du^\nu$ for $0\leq u \leq\sigma $.}

\medskip
It is obvious that the Beverton-Holt function  (\ref{a02}) satisfies the assumptions (B1) and (B2).

We claim here that problem (\ref{a11}) admits a global classical solution $v(t,y)$, that is, $v(t,y)$ is once continuously differentiable in $t$ ($t\in(0,+\infty)$), and twice continuously differentiable in $y$ ($y\in(0,l_0)$). In fact,  the initial value $v_0(y)\in C^1([0,l_0])$ and the fact that $g$ is once continuously differentiable implies that $v(0^+,y)\in C^1([0,l_0])$. By virtue of standard theory for parabolic equations, we can deduce that $v(t,y)\in C^{1,2}((0,T]\times(0,l_0))$.
Then, $v(T^+,y)=g(v(T,y))$ is also once continuously differentiable in $y$. Hence, let $v(T^+,y)$ be a new initial value for $t\in(T^+,2T]$, then $v(t,y)\in C^{1,2}((T,2T]\times(0,l_0))$. Eventually, we can obtain the solution $v(t,y)$ of problem (\ref{a11}) for $t \geq0$ and $y\in(0,l_0)$ by the same procedures.

Then we define
\[
PC([0,+\infty)\times [0,l_0])=\{v(t,y)\mid v(t,y)\in C((nT,(n+1)T]\times[0,l_0])\},
\]
\[
PC^{1,2}((0,+\infty)\times (0,l_0))=\{v(t,y)\mid v(t,y)\in C^{1,2}((nT,(n+1)T]\times(0,l_0))\}.
\]
The definition of upper and lower solutions  and the comparison principle for the initial boundary problem (\ref{a11}) with impulses are given as follows:
\begin{defi}
We say that $\tilde{v}(t,y),\hat{v}(t,y)\in PC^{1,2}((0,+\infty)\times(0,l_0))\cap PC([0,+\infty)\times[0,l_0])$ satisfying $0\leq \hat{v}(t,y)\leq \tilde{v}(t,y)$ are upper and lower solutions of problem $(\ref{a11})$, respectively, if $\tilde{v}(t,y)$ and $\hat{v}(t,y)$ make the following relationship true:
\small{\begin{eqnarray}
\left\{
\begin{array}{lll}
\tilde{v}_t\geq\frac{d}{\rho^2(t)}\tilde{v}_{yy}+(\alpha-\frac{\dot{\rho}(t)}{\rho(t)}) \tilde{v}-\gamma \tilde{v}^2,\; &\,y\in(0,l_0),\ t\in((nT)^+,(n+1)T], \\[2mm]
\hat{v}_t\leq\frac{d}{\rho^2(t)}\hat{v}_{yy}+(\alpha-\frac{\dot{\rho}(t)}{\rho(t)}) \tilde{v}-\gamma \hat{v}^2,\; &\,y\in(0,l_0),\ t\in((nT)^+,(n+1)T], \\[2mm]
\hat{v}(t,0)=0\leq \tilde{v}(t,0),\  \hat{v}(t,l_0)=0\leq \tilde{v}(t,l_0),\; &\,  t>0,\\[2mm]
\tilde{v}((nT)^+,y)\geq g(\tilde{v}(nT,y)),  \; &y\in(0,l_0),\ n=0,1,2,...\\
\hat{v}((nT)^+,y)\leq g(\hat{v}(nT,y)), \; &y\in(0,l_0),\ n=0,1,2,...
\end{array} \right.
\label{c01}
\end{eqnarray}}
with the initial condition
\bes
\begin{array}{llllll}
0\leq \hat{v}(0,y)\leq v_0(y)\leq \tilde{v}(0,y),\ \ y\in[0,l_0].
\end{array}
\label{c02}
\ees
\end{defi}

\begin{lem}(Comparison principle) Let $\tilde{v}(t,y)$ and $\hat{v}(t,y)$ be the upper and lower solutions to problem \eqref{a11},
then any solution $v(t,y)$ of problem $(\ref{a11})$ satisfies
\[
\hat{v}(t,y) \leq v(t,y)\leq \tilde{v}(t,y), \ t\in [0,\infty),\ y\in[0,l_0].
\]
\end{lem}
\textbf{Proof:} Let $\overline{w}(t,y)=\tilde{v}(t,y)-v(t,y)$. By the definition of an upper solution to problem (\ref{a11}) and assumption (B1), we have that:
\small{\[
\left\{
\begin{array}{lll}
\overline{w}(t,y)\geq \frac{d}{\rho^2(t)}\overline{w}_{yy}(t,y)+(\alpha-\frac{\dot{\rho}(t)}{\rho(t)}) \overline{w}(t,y)-\gamma (\overline{w}(t,y))^2, \; &\,y\in(0,l_0),\ t\in(0^+,T], \\[2mm]
\overline{w}(t,0) \geq 0,\ \overline{w}(t,l_0) \geq 0,\; &\,  t>0,\\[2mm]
\overline{w}(0^+,y)=\tilde{v}(0^+,y)-v(0^+,y)\geq g(\tilde{v}(0,y))-g(v_0(y))\geq 0,\; &\,y\in(0,l_0).

\end{array} \right.
\]}
The maximum principle implies that $\overline{w}(t,y)\geq 0$ for $t\in(0^+,T]$ and $y\in [0,l_0]$. Then, $\tilde{v}(t,y)\geq v(t,y)$ holds for $t\in(0^+,T]$
and $y\in [0,l_0]$. We take $\overline{w}(t,y)$ at time $t=T^+$ as a new initial value for the period $(T^+,2T]$,
which satisfies $\overline{w}(T^+,y)=\tilde{v}(T^+,y)-v(T^+,y)\geq g(\tilde{v}(T,y))-g(v(T,y))\geq 0$.
By the same procedures, $\tilde{v}(t,y)\geq v(t,y)$ can also be obtained for $t\in(T^+,2T]$ and $y\in [0,l_0]$.
Step by step, $\tilde{v}(t,y)\geq v(t,y)$ holds for all $t\geq 0$ and $y\in [0,l_0]$.

Similarly, the result that $v(t,y)\geq \hat{v}(t,y)$ for all $t\geq 0$ and $y\in [0,l_0]$ can also be deduced.           \qed\hfill

\medskip

The abovementioned preliminaries allow us to investigate the asymptotic behaviours of the solution to problem (\ref{a11}).
\begin{thm}
If $R_0<1$, then the solution $v(t,y)$ of problem $(\ref{a11})$ satisfies $\lim \limits_{t\rightarrow\infty}v(t,y)=0$ uniformly for $y\in [0,l_0]$.
\end{thm}
\textbf{Proof:}  Let
\[
\overline{v}(t,y)=Ke^{\alpha(1-\frac{1}{R_0})t}\phi(t,y),
\]
where $\phi(t,y)(\leq1)$ is a positive normalized eigenfunction of problem (\ref{b04}) and $K$ is a positive constant to be chosen later. It is easy to verify that $\overline{v}(t,y)$ is a solution of the following linear problem
\small{\begin{eqnarray}
\left\{
\begin{array}{lll}
v_t=\frac{d}{\rho^2(t)}v_{yy}+(\alpha-\frac{\dot{\rho}(t)}{\rho(t)}) v,\; &\,y\in(0,l_0),\ t\in((nT)^+,(n+1)T],\ n=0,1,2,..., \\[2mm]
v(t,0)=v(t,l_0)=0,\; &\,  t>0,\\[2mm]
v(0,y)=K\phi(0,y),\; &\, y\in(0,l_0),\\[2mm]
v((nT)^+,y)=g'(0)v(nT,y),\; &\,y\in(0,l_0),\ n=0,1,2,...
\end{array} \right.
\label{c03}
\end{eqnarray}}
It follows from (A1) that
\[
\frac{g(v)}{v}\leq \lim \limits_{\varepsilon\rightarrow0}\frac{g(\varepsilon)-g(0)}{\varepsilon}=g'(0)
\]
holds for sufficiently small $\varepsilon>0$,
which implies that
\bes
\begin{array}{llllll}
\overline{v}((nT)^+,y)=g'(0)\overline{v}(nT,y)\geq g(\overline{v}(nT,y)).
\end{array}
\label{c04}
\ees
Meanwhile, for any initial function $v(0,y)$, a sufficiently large constant $K$ can be chosen such that $\overline{v}(0,y)\geq v(0,y) $. Since the reaction term in problem (\ref{c03}) is larger than that in problem (\ref{a11}), we conclude that $\overline{v}(t,y)$ is a supersolution of problem (\ref{a11}). It then follows from Lemma 3.1 that
\[
v(t,y)\leq \overline{v}(t,y),\ \ t\geq 0,\ y\in[0,l_0].
\]
If $R_0<1$, then $\lim \limits_{t\rightarrow\infty}\overline{v}(t,y)=0$ for all $y\in[0,l_0]$. Therefore, $\lim \limits_{t\rightarrow\infty}v(t,y)=0$ uniformly for $y\in[0,l_0]$. \qed\hfill

\medskip
Next, in order to study the asymptotic behaviors of the solution to problem (\ref{a11}) when $R_0>1$, we study the following auxiliary periodic problem and state the relationship between its periodic solution and the solution of problem (\ref{a11}).
\begin{eqnarray}
\left\{
\begin{array}{lll}
v_t=\frac{d}{\rho^2(t)}v_{yy}+(\alpha-\frac{\dot{\rho}(t)}{\rho(t)}) v-\gamma v^2,\; &\,y\in(0,l_0),\ t\in((nT)^+,(n+1)T],\ n=0,1,2,..., \\[2mm]
v(t,0)=v(t,l_0)=0,\; &\,  t>0,\\[2mm]
v(0,y)=v(T,y),\; &\, y\in(0,l_0),\\[2mm]
v((nT)^+,y)=g(v(nT,y)),\; &\,y\in(0,l_0),\ n=0,1,2,...
\end{array} \right.
\label{c05}
\end{eqnarray}
The existence of a  solution to the periodic problem (\ref{c05}) is given by  upper and lower solutions technique.
\begin{defi}
We call $\tilde{v}(t,y),\hat{v}(t,y)\in PC^{1,2}((0,+\infty)\times(0,l_0))\cap PC([0,+\infty)\times[0,l_0])$ satisfying $0\leq \hat{v}(t,y)\leq \tilde{v}(t,y)$  an upper and lower solution of problem $(\ref{c05})$, respectively, if $\tilde{v}(t,y)$ and $\hat{v}(t,y)$ satisfy  relationship $(\ref{c01})$ and periodic conditions
\bes
\begin{array}{llllll}
\hat{v}(0,y)\leq \hat{v}(T,y),\ \tilde{v}(0,y)\geq \tilde{v}(T,y),\ y\in[0,l_0].
\end{array}
\label{c06}
\ees
\end{defi}

Obviously the upper and lower solutions of the periodic problem (\ref{c05}) are also the upper and lower solutions of the initial value problem (\ref{a11}) provided that $\hat{v}(t,y)\leq v(0,y) \leq \tilde{v}(t,y),y\in[0,l_0]$.

Let $f(v,t)=(\alpha-\frac{\dot{\rho}(t)}{\rho(t)})v-\gamma v^2$  and choose
\[
K^*=\frac{\alpha}{\gamma}+\sup_{t\in[0,T]}\frac{|\dot{\rho}(t)|}{\rho(t)}
\]
such that $F(v,t)=K^*v+f(v,t)$ is monotonically nondecreasing with respect to $v$. If there exists upper and lower solutions $\tilde{v}$ and $\hat{v}$ of problem (\ref{c05}), using $\overline{v}^{(0)}=\tilde{v}$ and $\underline{v}^{(0)}=\hat{v}$ as initial iteration, we can construct the iteration sequences $\{\overline{v}^{(m)}\}$ and $\{\underline{v}^{(m)}\}$ by the following process
\small{\begin{eqnarray}
\left\{
\begin{array}{lll}
\overline{v}^{(m)}_t-\frac{d}{\rho^2(t)}\overline{v}^{(m)}_{yy}+K^*\overline{v}^{(m)}= +K^*\overline{v}^{(m-1)} &\\
\qquad+(\alpha-\frac{\dot{\rho}(t)}{\rho(t)}) \overline{v}^{(m-1)}-\gamma (\overline{v}^{(m-1)})^2,\; &\,y\in(0,l_0),\ t\in((nT)^+,(n+1)T], \\[2mm]
\underline{v}^{(m)}_t-\frac{d}{\rho^2(t)}\underline{v}^{(m)}_{yy}+K^*\underline{v}^{(m)}=K^*\underline{v}^{(m-1)} &\\
\qquad +(\alpha-\frac{\dot{\rho}(t)}{\rho(t)}) \underline{v}^{(m-1)}-\gamma (\underline{v}^{(m-1)})^2,\; &\,y\in(0,l_0),\ t\in((nT)^+,(n+1)T], \\[2mm]
\overline{v}^{(m)}(t,0)=\overline{v}^{(m)}(t,l_0)=\underline{v}^{(m)}(t,0)=\underline{v}^{(m)}(t,l_0)=0,\; &\,  t>0,\\[2mm]
\overline{v}^{(m)}(0,y)=\overline{v}^{(m-1)}(T,y),\, \underline{v}^{(m)}(0,y)=\underline{v}^{(m-1)}(T,y),\; &\, y\in(0,l_0),\\[2mm]
\overline{v}^{(m)}((nT)^+,y)=g(\overline{v}^{(m-1)}((n+1)T,y)),\; &\,y\in(0,l_0),\ n=0,1,2,... ,\\[2mm]
\underline{v}^{(m)}((nT)^+,y)=g(\underline{v}^{(m-1)}((n+1)T,y)),\; &\,y\in(0,l_0),\ n=0,1,2,...
\end{array} \right.
\label{c07}
\end{eqnarray}}
Similarly as in \cite{pcv}, the monotone property of the iteration sequences is obtained in the following lemma.
\begin{lem}
If there exists upper and lower solutions $\tilde{v}$ and $\hat{v}$ of problem $(\ref{c05})$, then the sequences $\{\overline{v}^{(m)}\}$ and $\{\underline{v}^{(m)}\}$ have the monotone property
\[
\hat{v}\leq \underline{v}^{(m)} \leq \underline{v}^{(m+1)} \leq \overline{v}^{(m+1)} \leq \overline{v}^{(m)}\leq \tilde{v}
\]
for $t\in[0,+\infty)$ and $y\in[0,l_0]$.
\end{lem}
\textbf{Proof:}
Let $\underline{w}^{(0)}=\underline{v}^{(1)}-\underline{v}^{(0)}=\underline{v}^{(1)}-\hat{v}$.
By (\ref{c01}), (\ref{c06}) and (\ref{c07}), we have
\bes
\underline{w}^{(0)}(0,y)=\underline{v}^{(1)}(0,y)-\underline{v}^{(0)}(0,y)=\overline{v}^{(0)}(T,y)-\hat{v}(0,y)\geq0,\ y\in[0,l_0]
\label{c18}
\ees
and
\bes
\left\{
\begin{array}{lll}
\underline{w}^{(0)}-\frac{d}{\rho^2(t)}\underline{w}^{(0)}_{yy}+K^*\underline{w}^{(0)}=K^*\underline{v}^{(0)}+f(\underline{v}^{(0)},t)-(\hat{v}_t-
\frac{d}{\rho^2(t)}\hat{v}_{yy}+K^*\hat{v})\\
\ \ \ \ \ \ \ \ \ \ \ \ \ \ \ \ \ \ \ \ \ \ \ \ \ \ \ \ \ \ \ \ \ \ \ \ \ \ =f(\hat{v},t)-(\hat{v}_t-
\frac{d}{\rho^2(t)}\hat{v}_{yy})\geq0,\quad y\in(0,l_0),\ t\in(0^+,T], \\[2mm]
\underline{w}^{(0)}(t,0)=\underline{v}^{(1)}(t,0)-\hat{v}(t,0)=0,\ \underline{w}^{(0)}(t,l_0)=\underline{v}^{(1)}(t,l_0)-\hat{v}(t,l_0)=0,\quad t>0,\\[2mm]
\underline{w}^{(0)}(0^+,y)=\underline{v}^{(1)}(0^+,y)-\hat{v}(0^+,y)=g(\underline{v}^{(0)}(T,y))-g(\hat{v}(0,y))\geq0,\quad y\in(0,l_0).
\end{array} \right.
\label{c19}
\ees
(\ref{c19}) and the positivity lemma for parabolic problems show that $\underline{w}^{(0)}\geq0$ for $t\in(0,T]$ and $y\in[0,l_0]$, which together with (\ref{c18}) yields that $\underline{w}^{(0)}\geq0$ for $t\in[0,T]$ and $y\in[0,l_0]$ i.e. $\underline{v}^{(1)}\geq \underline{v}^{(0)}$ for $t\in[0,T]$ and $y\in[0,l_0]$. Then, if $\underline{w}^{(0)}(T^+,y)$ is a new initial value for $t\in(T,2T]$, we can deduce that $\underline{v}^{(1)}\geq \underline{v}^{(0)}$ for $t\in(T,2T]$ and $y\in[0,l_0]$ similarly. Hence $\underline{v}^{(1)}\geq \underline{v}^{(0)}$ holds for $t\in[0,+\infty)$ and $y\in[0,l_0]$. Similarly, the property of an upper solution gives the result that $\overline{v}^{(1)}\geq \overline{v}^{(0)}$ for $t\in[0,+\infty)$ and $y\in[0,l_0]$. Moreover, let
$w^{(1)}=\overline{v}^{(1)}-\underline{v}^{(1)}$, then
\[
\left\{
\begin{array}{lll}
w^{(1)}-\frac{d}{\rho^2(t)}w^{(1)}_{yy}+K^*w^{(1)}=F(\overline{v}^{(0)},t)-F(\underline{v}^{(0)},t)\geq0,\quad y\in(0,l_0),\ t\in(0^+,T], \\[2mm]
w^{(1)}(t,0)=\overline{v}^{(1)}(t,0)-\underline{v}^{(1)}(t,0)=0,\ w^{(1)}(t,l_0)=\overline{v}^{(1)}(t,l_0)-\underline{v}^{(1)}(t,l_0)=0,\quad t>0,\\[2mm]
w^{(1)}(0,y)=\overline{v}^{(1)}(0,y)-\underline{v}^{(1)}(0,y)=\tilde{v}(T,y)-\hat{v}(T,y)\geq0,\quad y\in(0,l_0),\\[2mm]
w^{(1)}(0^+,y)=\overline{v}^{(1)}(0^+,y)-\underline{v}^{(1)}(0^+,y)=g(\overline{v}^{(0)}(T,y))-g(\underline{v}^{(0)}(T,y))\geq0,\quad y\in(0,l_0),
\end{array} \right.
\]
which yields that $w^{(1)}\geq0$ for $t\in[0,T]$ and $y\in[0,l_0]$. The result that $w^{(1)}\geq0$ for $t\in[0,+\infty)$ and $y\in[0,l_0]$ can be deduced by the same procedures. Therefore, we obtain that $\underline{v}^{(0)}\leq \underline{v}^{(1)} \leq \overline{v}^{(1)} \leq \overline{v}^{(0)}$ for $t\in[0,+\infty)$ and $y\in[0,l_0]$. By the method of induction, the sequences have the following property
\[
\hat{v}\leq \underline{v}^{(m)} \leq \underline{v}^{(m+1)} \leq \overline{v}^{(m+1)} \leq \overline{v}^{(m)}\leq \tilde{v}
\]
for $t\in[0,+\infty)$ and $y\in[0,l_0]$. \qed\hfill

Now we investigate the existence and uniqueness of a positive periodic solution of problem (\ref{c05}).

\begin{thm}
If $R_0>1$, then problem $(\ref{c05})$ admits a unique positive periodic solution $v^*(t,y)$.
\end{thm}
\textbf{Proof:}
First of all, we construct the upper solution of problem (\ref{c05}).
Let $\tilde{v}=MW(t)\ (M>1)$ with $W(t)$ satisfying
\begin{eqnarray}
\left\{
\begin{array}{lll}
W_t(t)=(\alpha-\frac{\dot{\rho}(t)}{\rho(t)})W(t)-\gamma W^2(t), \; &\, t\in((nT)^+,(n+1)T], \ n=0,1,2,... ,\\[2mm]
W(t)=W(t+T),\; &\,t\geq 0,\\[2mm]
W((nT)^+)=g'(0)W(nT)\geq g(W(nT)),\;&\, n=0,1,2,...
\end{array} \right.
\label{c08}
\end{eqnarray}
It is clear that $\tilde{v}=MW(t)\ (M>1)$ is an upper solution of problem (\ref{c05}).
Problem (\ref{c08}) can be solved by direct calculations.
Integrating from $(nT)^+$ to $t (t\in((nT)^+,(n+1)T])$ in the first equation in (\ref{c08}) yields
\bes
\begin{array}{llllll}
W(t)=\frac{e^{\alpha t}W((nT)^+)}{W((nT)^+)\int_{nT}^t \frac{\gamma e^{\alpha\tau}}{\rho(\tau)}d\tau+e^{\alpha nT}}, \ t\in ((nT)^+,(n+1)T],
\end{array}
\label{c09}
\ees
then
\[
\begin{array}{llllll}
W((n+1)T)&=\frac{e^{\alpha (n+1)T}g'(0)W(nT)}{g'(0)W(nT)\int_{nT}^{(n+1)T} \frac{\gamma e^{\alpha\tau}}{\rho(\tau)}d\tau+e^{\alpha nT}}
=\frac{e^{\alpha (n+1)T}g'(0)W(nT)}{g'(0)W(nT)\int_{0}^{T} \frac{\gamma e^{\alpha(\tau+nT)}}{\rho(\tau)}d\tau+e^{\alpha nT}}\\
&=\frac{e^{\alpha T}g'(0)W(nT)}{g'(0)W(nT)\int_{0}^{T} \frac{\gamma e^{\alpha\tau}}{\rho(\tau)}d\tau+1},
\end{array}
\]
Since $R_0>1$, then $\alpha T>d\lambda_1\int_0^T\frac{1}{\rho^2(t)}dt-\ln g'(0)>-\ln g'(0)$, that is, $e^{\alpha T}g'(0)>1$.
Using the periodicity, we obtain that
 $W(nT)= \frac{e^{\alpha T}g'(0)-1}{g'(0)\int_0^T \frac{ \gamma e^{\alpha\tau}}{\rho(\tau)}d\tau}>0$ and
\[
W(t)=\frac{e^{\alpha t}(e^{\alpha T}g'(0)-1)}{(e^{\alpha T}g'(0)-1)\int_{nT}^{t} \frac{\gamma e^{\alpha\tau}}{\rho(\tau)}d\tau+e^{\alpha nT}\int_{0}^{T} \frac{\gamma e^{\alpha\tau}}{\rho(\tau)}d\tau},\ t\in ((nT)^+,(n+1)T].
\]
So now we have the upper solution $\tilde{v}$ of problem (\ref{c05}).

Next we consider the lower solution and define
\begin{eqnarray}
\hat{v}(t,y)=\left\{
\begin{array}{lll}
\varepsilon \phi(nT,y), \; &\, t=nT, \\[2mm]
\varepsilon\frac{\rho_1}{g'(0)}\phi((nT)^+,y),\; &\,t=(nT)^+,\\[2mm]
\varepsilon\frac{\rho_1}{g'(0)} e^{[\alpha(1-\frac{1}{R_0})-\delta](t-nT)}\phi(t,y),\;&\,t\in((nT)^+,(n+1)T],\  n=0,1,2,... ,
\end{array} \right.
\label{c10}
\end{eqnarray}
where the positive eigenfunction $\phi(t,y)$ is defined in (\ref{b04}) and $\varepsilon(>0)$ is a sufficiently small constant to be chosen later, as well as positive constants $\delta=\frac{\alpha}{2}(1-\frac{1}{R_0})$ and $\rho_1=e^{-\frac{\alpha}{2}(1-\frac{1}{R_0})T}g'(0)$ are chosen to make sure that $\hat{v}(nT,y)=\hat{v}((n+1)T,y)$.
For $t \in ((nT)^+,(n+1)T]$ and $y\in (0,l_0)$, if $\varepsilon<\varepsilon_1:=\frac{\delta}{\gamma}$, we have
\[
\begin{array}{llllll}
\frac{\partial \hat{v}}{\partial t}-[\frac{d}{\rho^2(t)}\hat{v}+(\alpha-\frac{\dot{\rho}(t)}{\rho(t)})\hat{v}-\gamma \hat{v}^2]
\\
=[\alpha(1-\frac{1}{R_0})-\delta]\varepsilon\frac{\rho_1}{g'(0)}e^{[\alpha(1-\frac{1}{R_0})-\delta](t-nT)}\phi+\varepsilon\frac{\rho_1}{g'(0)} e^{[\alpha(1-\frac{1}{R_0})-\delta](t-nT)}[\frac{d}{\rho^2(t)}\phi_{yy}+(\frac{\alpha}{R_0}-\frac{\dot{\rho(t)}}{\rho(t)})\phi]
\\
\quad -[\frac{d}{\rho^2(t)}\varepsilon\frac{\rho_1}{g'(0)} e^{[\alpha(1-\frac{1}{R_0})-\delta](t-nT)}\phi_{yy}+(\alpha-\frac{\dot{\rho}(t)}{\rho(t)})\varepsilon\frac{\rho_1}{g'(0)} e^{[\alpha(1-\frac{1}{R_0})-\delta](t-nT)}\phi]\\
\quad +\gamma(\varepsilon\frac{\rho_1}{g'(0)} e^{[\alpha(1-\frac{1}{R_0})-\delta](t-nT)}\phi)^2
\\
=\varepsilon[-\delta+\gamma\varepsilon\frac{\rho_1}{g'(0)} e^{[\alpha(1-\frac{1}{R_0})-\delta](t-nT)}\phi]\frac{\rho_1}{g'(0)} e^{[\alpha(1-\frac{1}{R_0})-\delta](t-nT)}\phi
\\
<0.
\end{array}
\]
On the other hand it follows from assumption (B2) that
\[
\begin{array}{llllll}
g(\hat{v}(nT,y))-\hat{v}((nT)^+,y)\\
=g(\hat{v}(nT,y))-\varepsilon\frac{\rho_1}{g'(0)}\phi((nT)^+,y)

=g(\hat{v}(nT,y))-\rho_1\hat{v}(nT,y)\\
\geq (g'(0)-\rho_1)\hat{v}(nT,y)-D\{\varepsilon\phi(nT,y)  \}^\nu
\\
=[(g'(0)-\rho_1)-D(\varepsilon\phi(nT,y))^{\nu-1}]\varepsilon\phi(nT,y) \\
\geq 0
\end{array}
\]
if $\varepsilon<\varepsilon_2:=(\frac{g'(0)-\rho_1}{D})^\frac{1}{\nu-1}$.
These results show that $\hat{v}(t,y)$ is a lower solution of problem (\ref{c05}).

Next, using $\overline{v}^{(0)}=\tilde{v}$ and $\underline{v}^{(0)}=\hat{v}$ as the initial iteration, two sequences $\{ \overline{v}^{(m)}\}$ and $\{\underline{v}^{(m)}\} $ are constructed from problem (\ref{c07}). By Lemma 3.2, we have that
\[
\hat{v}\leq \underline{v}^{(m)} \leq \underline{v}^{(m+1)} \leq \overline{v}^{(m+1)} \leq \overline{v}^{(m)}\leq \tilde{v}
\]
for every $m=1,2,...$  Hence, the limits of sequences $\{ \overline{v}^{(m)}\}$ and $\{\underline{v}^{(m)}\} $ exist and
\[
\lim_{m\rightarrow\infty}\overline{v}^{(m)}=\overline{v}^*,\ \lim_{m\rightarrow\infty}\underline{v}^{(m)}=\underline{v}^*,
\]
where $\overline{v}^*$ and $\underline{v}^*$ are T-periodic solutions of problem (\ref{c05}).
Furthermore,
\[
\hat{v}\leq \underline{v}^{(m)} \leq \underline{v}^{(m+1)} \leq \underline{v}^* \leq \overline{v}^* \leq \underline{v}^{(m+1)} \leq \underline{v}^{(m)}\leq \tilde{v}.
\]
Now we claim that $\overline{v}^*$ and $\underline{v}^*$ are the maximal and minimal positive periodic solutions of problem (\ref{c05}).  For any positive periodic solution $v(t,y)$ of problem (\ref{c05}) satisfying $\hat{v}\leq v \leq \tilde{v}$, we use the same iteration as in (\ref{c07}) with the initial iteration $\overline{v}^{(0)}=\tilde{v}$ and $\underline{v}^{(0)}=v$, where $\tilde{v}$ and $v$ are a pair of upper and lower solutions of problem (\ref{c05}), from which we can deduce that
\[
v(t,y)\leq \bar{v}^*(t,y),\ t\geq0,\ y\in[0,l_0],
\]
so that  $\overline{v}^*$ is the maximal positive periodic solution of problem (\ref{c05}). Similarly, $\underline{v}^*$ is the minimal positive periodic solution of problem (\ref{c05}).

Now the proof ends with showing the uniqueness of the positive periodic solution of problem (\ref{c05}).
Let $v_1$ and $v_2$ be two solutions, and define sector
\[
\Lambda=\{ s\in[0,1], sv_1 \leq v_2,\  t=0,\ t=0^+,\ t\in(0^+,T],\ y\in[0,l_0] \}
\]
It is easily seen that $\Lambda$ contains a neighbourhood near by $0$. We claim that $1\in \Lambda$. If not, we assume that $s_0=\sup \Lambda<1$.
Recalling $F(v,t)=f(v,t)+K^*v$ is nondecreasing and $\frac{f(v,t)}{v}$ is decreasing in $v$ on $[0,\max v_2]$, we deduce that
\[
\begin{array}{llllll}
(v_2-s_0v_1)_t-\frac{d}{\rho^2(t)}(v_2-s_0v_1)_{yy}+K^*(v_2-s_0v_1)\\
\ =f(v_2,t)+K^*v_2-s_0(f(v_1,t)+K^*v_1)\\
\ \geq f(s_0v_1,t)+K^*s_0v_1-s_0(f(v_1,t)+K^*v_1)\geq0
\end{array}
\]
for $t\in(0^+,T]$ and $y\in(0,l_0)$. Using assumptions (A1) and (B1) we find that
\[
\begin{array}{llllll}
v_2(0^+,y)-s_0v_1(0^+,y)=g(v_2(0,y))-s_0g(v_1(0,y))\\
\ \geq g(s_0v_1(0,y))-s_0g(v_1(0,y))\geq 0
\end{array}
\]
for $y\in(0,l_0)$.
On the other hand, for $t>0$,
\[
v_2(t,0)-s_0v_1(t,0)=v_2(t,l_0)-s_0v_1(t,l_0)=0.
\]
By the strong maximum principle \cite{pmh}, we have the following assertions:

(i) $v_2-s_0v_1>0$ holds for $t=0^+$, $t\in(0^+,T]$ and $y\in(0,l_0)$. Recalling $v_1$ and $v_2$ are T-periodic solutions, that is $v_1(0,y)=v_1(T,y)$ and $v_2(0,y)=v_2(T,y)$ for $y\in(0,l_0)$, then $v_2-s_0v_1>0$ holds for  $t\in[0,T]$ and $y\in(0,l_0)$. By Hopf's boundary lemma, $\frac{\partial}{\partial\eta}\mid_{y=0}(v_2-s_0v_1)>0$ and $\frac{\partial}{\partial\eta}\mid_{y=l_0}(v_2-s_0v_1)<0$ for $t\in[0,T]$, where $\eta$ is the outward unit normal vector. Then there exists a constant $\epsilon>0$ such that $v_2-s_0v_1\geq\epsilon v_1$, which leads to $s_0+\epsilon\in\Lambda$. This contradicts  the maximum property of $s_0$.

(ii) $v_2-s_0v_1\equiv0$ holds for $t=0^+$, $t\in(0^+,T]$ and $y\in(0,l_0)$. We then obtain that $f(v_2,t)=s_0f(v_1,t)$. However, since $s_0<1$, $f(v_2,t)=f(s_0v_1,t)>s_0f(v_1,t)$. Hence this case is impossible.

We conclude that problem (\ref{c05}) has a unique positive periodic solution $v^*(t,y)$.\qed  \hfill

\medskip
We have now established the existence and uniqueness of a positive periodic solution of problem (\ref{c05})  in the above theorem. Now we show how this periodic solution is a global attractor for the problem (\ref{a11}) This shown in the following theorem .
\begin{thm}
If $R_0>1$, then for nonnegative nontrivial initial value $v_0(y)$,
any solution $v(t,y)$ of problem $(\ref{a11})$  satisfies
\[
\lim_{m\rightarrow\infty}v(t+mT,y)\rightarrow v^*(t,y),\ t\geq0,\ y\in[0,l_0],
\]
where $v^*(t,y)$ is a positive T-periodic solution of problem $(\ref{c05})$.
That is, $v^*(t,y)$ is a global attractor of problem $(\ref{a11})$.
\end{thm}
\textbf{Proof:}
Without loss of generality, we assume that $v_0(y)>0$ for $y\in(0,l_0)$. Otherwise we can replace the initial time 0 by any $t_0>0$ since $v_0(t_0,y)>0$ for $y\in(0,l_0)$. Noting that $\phi_y(0,0)>0$ and $\phi_y(0,l_0)<0$, by Hopf's boundary lemma, we can choose a sufficiently small $\varepsilon$ such that $\varepsilon\phi(0,y)\leq v(0,y)$. Also, a sufficiently big $M$ can be chosen such that $v(0,y) \leq MW(0)$. For given $\varepsilon$ and $M$, the function $\tilde{v}:=MW(t)$ with $W(t)$ defined in (\ref{c08}) and $\hat{v}$ defined in (\ref{c10}), satisfies
\bes
\begin{array}{llllll}
\hat{v}(0,y)\leq v(0,y) \leq \tilde{v}(0,y),\ y\in[0,l_0].
\end{array}
\label{c11}
\ees
Since $g$ is nondecreasing with respect to $v$, we have that
\[
\hat{v}(0^+,y)=g(\hat{v}(0,y))\leq g(v(0,y))=v(0^+,y)\leq g(\tilde{v}(0,y))=\tilde{v}(0^+,y).
\]
It follows from the classical comparison principle that $\hat{v}(t,y)\leq v(t,y) \leq \tilde{v}(t,y), \ t\in (0^+,T], y\in[0,l_0]$.
Induction shows that $\hat{v}(t,y)\leq v(t,y) \leq \tilde{v}(t,y), \ \ t=nT,\ (nT)^+,\ t\in ((nT)^+,(n+1)T], \ y\in[0,l_0]$. Hence, for $n=0,1,2,...$,
 \bes
\underline{v}^{(0)}(t,y)\leq v(t,y) \leq \overline{v}^{(0)}(t,y), \ \ t=nT,\ (nT)^+,\ t\in ((nT)^+,(n+1)T], \ y\in[0,l_0].
\label{c12}
\ees
Then,
\bes
\underline{v}^{(0)}(T,y)\leq v(T,y) \leq \overline{v}^{(0)}(T,y),\ y\in[0,l_0],
\label{c21}
\ees
which together with $\underline{v}^{(1)}(0,y)=\underline{v}^{(0)}(T,y)$ and $\overline{v}^{(1)}(0,y)=\overline{v}^{(0)}(T,y)$ yields
\[
\underline{v}^{(1)}(0,y)\leq v(T,y) \leq \overline{v}^{(1)}(0,y),\ y\in[0,l_0].
\]
By the monotonicity of $g$ and (\ref{c21}), we obtain
\[
g(\underline{v}^{(0)}(T,y))\leq g(v(T,y)) \leq g(\overline{v}^{(0)}(T,y)),\ y\in[0,l_0].
\]
Due to the last two equations in (\ref{c07}) and the last equation in (\ref{a11}) we have that
\[
\underline{v}^{(1)}(0^+,y)=g(\underline{v}^{(0)}(T,y))\leq g(v(T,y))=v(T^+,y) \leq g(\overline{v}^{(0)}(T,y))=\overline{v}^{(1)}(0^+,y),\ y\in[0,l_0],
\]
that is,
\[
\underline{v}^{(1)}(0^+,y)\leq v(T^+,y)\leq \overline{v}^{(1)}(0^+,y),\ y\in[0,l_0].
\]
Then,  from the classical comparison principle  $\underline{v}^{(1)}(t,y)\leq v(t+T,y)\leq \overline{v}^{(1)}(t,y)$ holds for $t\in (0^+,T]$ and $y\in[0,l_0]$. By  induction we have
\[
\underline{v}^{(1)}(t,y)\leq v(t+T,y)\leq \overline{v}^{(1)}(t,y),\ \ t=nT,\ (nT)^+,\ t\in ((nT)^+,(n+1)T], \ y\in[0,l_0]
\]
for $n=0, 1, 2, ...$ Also by induction, it follows from the last two equations in problem (\ref{c07}) that
\bes
\underline{v}^{(m)}(t,y)\leq v(t+mT,y)\leq \overline{v}^{(m)}(t,y),\ t\geq0, y\in[0,l_0],
\label{c22}
\ees
since (\ref{c22}) holds for $m=0$ and $m=1$.
With $\lim\limits_{m\rightarrow \infty} \underline{v}^{(m)}(t,y)=\lim\limits_{m\rightarrow \infty} \overline{v}^{(m)}(t,y)=v^*(t,y)$  by the uniqueness of the periodic solution of problem (\ref{c05}) presented in Theorem 3.4, we conclude that
\[
\lim_{m\rightarrow\infty}v(t+mT,y)\rightarrow v^*(t,y),\ t\geq0,\ y\in[0,l_0].    \qed  \hfill
\]

\subsection{Nonmonotone Case}
To investigate the nonmonotone case we need the following restriction, just  like in \cite{lma}:

\medskip
(B3) \ There exists a $\sigma_0>0$ such that $g(u)$ is nondecreasing for $0\leq u \leq \sigma_0$.

\medskip
This assumption is consistent with the actual harvesting mechanism.  The classical pluse function, e.g the Ricker function $g(u)=ue^{r-bu}$, satisfies this;  $g(u)$ is increasing for $0<u<1/b$ and decreasing for $u>1/b$.

We define
\bes
\begin{array}{llllll}
g^+(u)=\max\limits_{0\leq w \leq u}g(w)
\end{array}
\label{c13}
\ees
for $u\geq 0$. Obviously, for $u\geq0$, $g^+(u)$ is nondecreasing and $g^+(u)\geq g(u)$, $g^+(u)=g(u)$ for small $u>0$ and ${g^+}'(0)=g'(0)$. Under the condition $R_0>1$, the following problem admits a minimal positive periodic solution $\underline{v^+}(t)$ satisfying
\begin{eqnarray}
\left\{
\begin{array}{lll}
v_t(t)=(\alpha-\frac{\dot{\rho}(t)}{\rho(t)})v(t)-\gamma v^2(t), \; &\, t\in((nT)^+,(n+1)T], \ n=0,1,2,... ,\\[2mm]
v(0)=v(T),\; &\,t\geq 0,\\[2mm]
v((nT)^+)=g^+(v(nT)),\;&\, n=0,1,2,...
\end{array} \right.
\label{c20}
\end{eqnarray}
Let $\beta^+=\min\limits_{t\in[0,T]}\underline{v^+}(t)$.
We also define
\bes
\begin{array}{llllll}
g^-(u)=\min\limits_{u\leq w \leq \beta^+}g(w)
\end{array}
\label{c14}
\ees
for $0\leq u\leq \beta^+$. It is obvious to see that for $u\geq0$, $g^-(u)$ is nondecreasing and $g^-(u)\leq g(u)$, $g^-(u)=g(u)$ for small $u>0$ and ${g^-}'(0)=g'(0)$.

Then we take the following two auxiliary problems into account
\begin{eqnarray}
\left\{
\begin{array}{lll}
v_t=\frac{d}{\rho^2(t)}v_{yy}+(\alpha-\frac{\dot{\rho}(t)}{\rho(t)}) v-\gamma v^2,\; &\,y\in(0,l_0),\ t\in((nT)^+,(n+1)T],\ n=0,1,2,... \\[2mm]
v(t,0)=v(t,l_0)=0,\; &\,  t>0,\\[2mm]
v(0,y)=v_0(y)\geq,\not\equiv0,\; &\, y\in[0,l_0],\\[2mm]
v((nT)^+,y)=g^+(v^+(nT,y)),\; &\,y\in(0,l_0),\ n=0,1,2,... ,
\end{array} \right.
\label{c15}
\end{eqnarray}
and
\begin{eqnarray}
\left\{
\begin{array}{lll}
v_t=\frac{d}{\rho^2(t)}v_{yy}+(\alpha-\frac{\dot{\rho}(t)}{\rho(t)}) v-\gamma v^2,\; &\,y\in(0,l_0),\ t\in((nT)^+,(n+1)T],\ n=0,1,2,... \\[2mm]
v(t,0)=v(t,l_0)=0,\; &\,  t>0,\\[2mm]
v(0,y)=v_0(y)\geq,\not\equiv0,\; &\, y\in[0,l_0],\\[2mm]
v((nT)^+,y)=g^-(v^-(nT,y)),\; &\,y\in(0,l_0),\ n=0,1,2,... ,
\end{array} \right.
\label{c16}
\end{eqnarray}
We use $v^+(t,y)$ and $v^-(t,y)$ to represent the solution of problem (\ref{c15}) and (\ref{c16}), respectively.
Recalling that $g^-(v)\leq g(v) \leq g^+(v)$ for $u\geq0$, Lemma 3.1 implies that if $0\leq v_0^+(y)\leq v_0(y)\leq v_0^-(y)\leq \beta^+$, then any solution $v(t,y)$ of problem (\ref{a11}) satisfies
\bes
\begin{array}{llllll}
0\leq v^-(t,y)\leq v(t,y)\leq v^+(t,y)\leq \beta^+.
\end{array}
\label{c17}
\ees
By expression (\ref{b05}), problem (\ref{c15}) and (\ref{c16}) have the same basic reproduction number $R_0$. By Theorems 3.2, 3.4 and 3.5,
 problem (\ref{c15}) and (\ref{c16}) have the same asymptotic behaviors of solutions. Owing to this and (\ref{c17}), when $R_0<1$, we obtain that the solution $v(t,y)$ of problem $(\ref{a11})$ converges to 0 as $t\rightarrow \infty$. If $R_0>1$, then the solution $v^-(t,y)$ of problem (\ref{c16})  satisfies
\[
\underline{v^-}(t,y)\leq \liminf_{m\rightarrow+\infty}v^-(t+mT,y)\leq \limsup_{m\rightarrow+\infty}v^-(t+mT,y)\leq \overline{v^-}(t,y),\ t\geq0,\ y\in[0,l_0],
\]
where $\underline{v^-}(t,y)$ and $\overline{v^-}(t,y)$ are the minimum and maximum positive periodic solutions of the corresponding periodic
problem related to problem (\ref{c16}), respectively.
By (\ref{c17}), we can deduce that any positive solution $v(t,y)$ of problem (\ref{a11}) satisfies
\[
\liminf_{m\rightarrow+\infty}v(t+mT,y)\geq\underline{v^-}(t,y).
\]
Hence, the asymptotic behavior of the solution to problem (\ref{a11}) can be given as follows:
\begin{thm}
$(i)$ If $R_0<1$, then the solution $v(t, y)$ of problem $(\ref{a11})$ satisfies $\lim \limits_{t\rightarrow\infty} v(t, y)=0$.\\
$(ii)$ If $R_0>1$, then any solution $v(t, y)$ of problem $(\ref{a11})$  satisfies
\[
\liminf_{m\rightarrow +\infty} v(t+mT, y)\geq \underline {v^-}(t,y)
\]
provided that $0 \not\equiv, \leq v_0(y) \leq \beta^+$, where $\underline{v^-}(t, y)$ is the minimal positive periodic solution of
the corresponding periodic problem related to problem (\ref{c16}).
\end{thm}

\section{\bf The effects of the evolution rate and impulsive harvesting}
The above sections show the threshold-type dynamics for a population,  determined by diffusion parameters, evolution rate and properties of
the impulsive harvesting. In this section, we aim to investigate how the evolution rate and impulsive harvesting affect these dynamical
behaviors of a population by adopting numerical analyses. In all simulations, we consider the interval $[0,l(t)]=[0,\rho(t)l_0]$,
where $l_0=\pi$, and we fix some parameters
\[
d=1, \ \ \alpha=1.1,\ \  \gamma=0.05
\]
and subsequently obtain $\lambda_1=(\frac{\pi}{l_0})^2=1$. We  choose $v_0(y)=0.5\sin(x)+0.2\sin(3x)$ as initial function.

\subsection{The effect of the evolution rate}
\begin{figure}[ht]
\centering
\subfigure[]{ {
\includegraphics[width=0.3\textwidth]{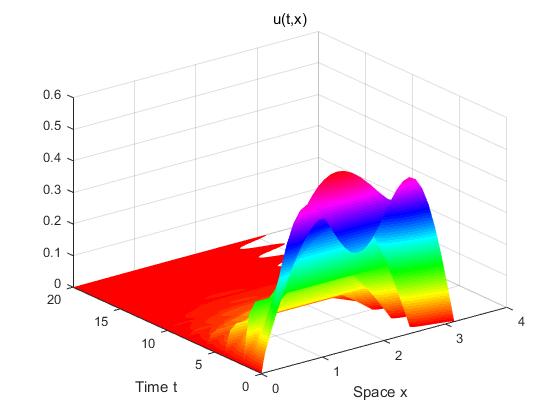}
} }
\subfigure[]{ {
\includegraphics[width=0.3\textwidth]{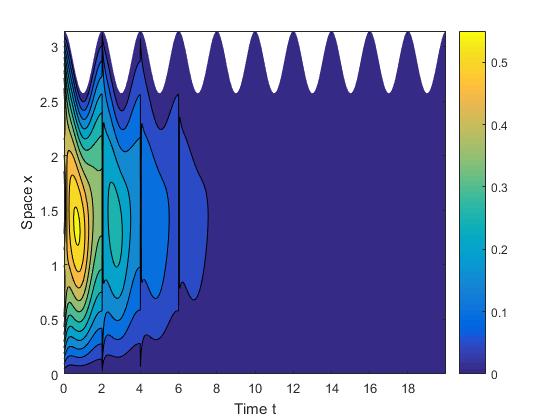}
} }
\subfigure[]{ {
\includegraphics[width=0.3\textwidth]{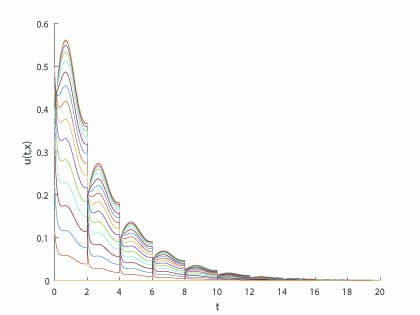}
} }
\caption{\scriptsize  $\rho_1(t)=e^{-0.1(1-\cos\pi t)}$ and $g(u)$ equals the Beverton-Holt function with $a=10$ and $m=8$.
The domain is periodically evolving with $\rho_1$ and $R_0<1$. Graphs $(a)-(c)$ show that the population $u(t,x)$ decays to $0$. Graphs $(b)$ and $(c)$ are the
cross-sectional view and projection of $u$ on the $t-u$-plane, respectively. The color bar in graph $(b)$ shows the density of $u(t,x)$.}
\label{tu1}
\end{figure}

\begin{figure}[ht]
\centering
\subfigure[]{ {
\includegraphics[width=0.30\textwidth]{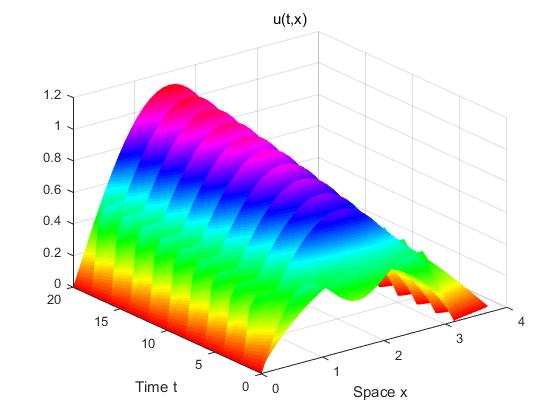}
} }
\subfigure[]{ {
\includegraphics[width=0.30\textwidth]{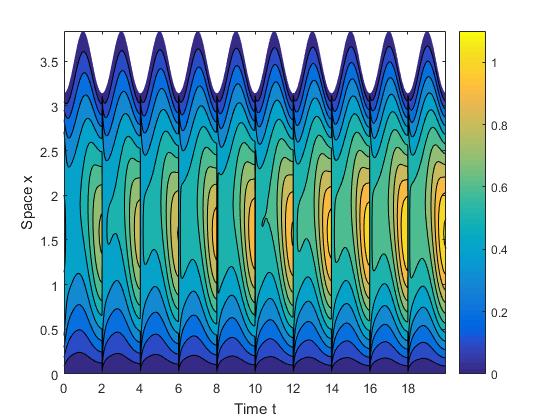}
} }
\subfigure[]{ {
\includegraphics[width=0.30\textwidth]{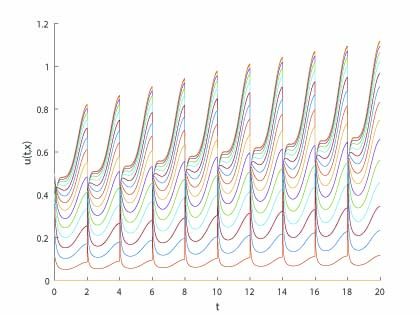}
} }
\caption{\scriptsize  $\rho_2(t)=e^{0.1(1-\cos\pi t)}$ and $g(u)$ equals the Beverton-Holt function with $a=10$ and $m=8$, then $R_0>1$.
Spatiotemporal dynamics of population $u(t,x)$ are shown in graph $(a)$,  indicating that the population approaches a positive periodic steady state. Graph $(b)$ is the
cross-sectional view and indicates the periodic evolution of the domain. The appearance of impulsive harvesting every time $T=2$
can be seen in graph $(c)$, which is the projection of $u$ on the $t-u$-plane. }
\label{tu2}
\end{figure}

We  choose two different $\rho(t)$ to analyze the effect of the evolution rate on the dynamics of the population when the impulsive
harvesting occurs every time $T=2$. We first fix $g(u)$ as the Beverton-Holt function with $a=10$ and $m=8$ as in (\ref{a02}),
in which case $g(u)$ is a monotonically increasing function.
\begin{exm}
Fix $g(u)=8u/(10+u)$ and then $g'(0)=0.8$. We first choose $\rho_1(t)=e^{-0.1(1-\cos\pi t)}$, it follows from \eqref{b05} that
$$R_0(\rho_1)=\frac{\alpha}{\frac{d\lambda_1}{T}\int_0^T\frac{1}{\rho^2(t)}dt-\frac{1}{T}\ln g'(0)}\approx0.8177<1.$$
One can see from Fig. \ref{tu1} that the population suffers eventual extinction.

Now we choose $\rho_2(t)=e^{0.1(1-\cos\pi t)}$. Direct calculations show that
$$R_0(\rho_2)=\frac{\alpha}{\frac{d\lambda_1}{T}\int_0^T\frac{1}{\rho^2(t)}dt-\frac{1}{T}\ln g'(0)}\approx1.1721>1.$$
It is easily seen from Fig. \ref{tu2} that the population approaches a positive periodic steady state.

The example shows that the population vanishes in a periodically evolving habitat with a small evolution rate and persists in a habitat with a larger evolution rate.
\end{exm}

\begin{figure}[ht]
\centering
\subfigure[]{ {
\includegraphics[width=0.30\textwidth]{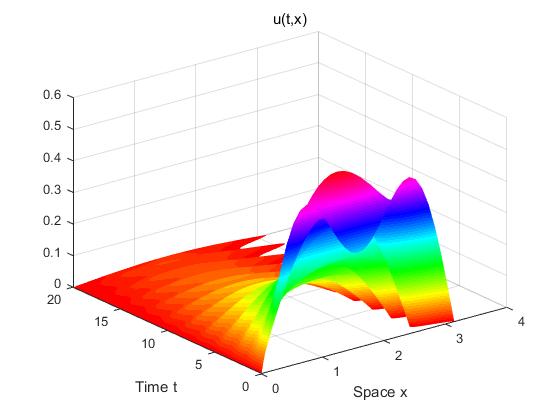}
} }
\subfigure[]{ {
\includegraphics[width=0.30\textwidth]{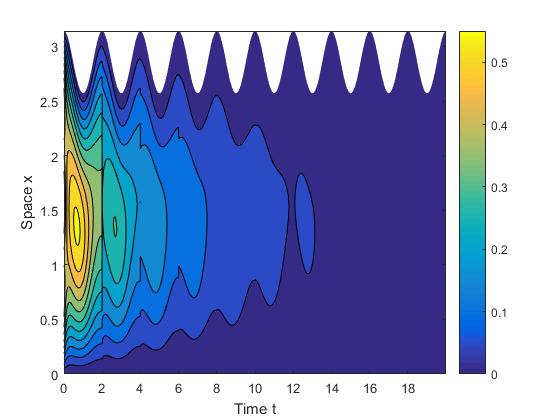}
} }
\subfigure[]{ {
\includegraphics[width=0.30\textwidth]{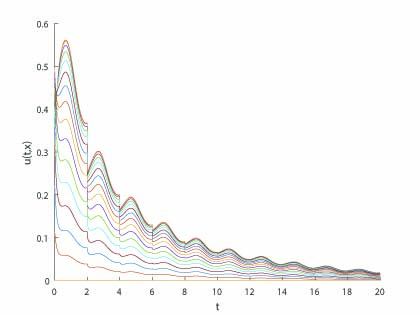}
} }
\caption{\scriptsize $\rho_1(t)=e^{-0.1(1-\cos\pi t)}$ and $g(u)$ equals the Ricker function with $b=1.2$ and $r=0.05$. The domain is
periodically evolving with $\rho_1$ and $R_0<1$. Graphs $(a)-(c)$ show that the population $u(t,x)$ decays to $0$. }
\label{tu3}
\end{figure}

\begin{figure}[ht]
\centering
\subfigure[]{ {
\includegraphics[width=0.30\textwidth]{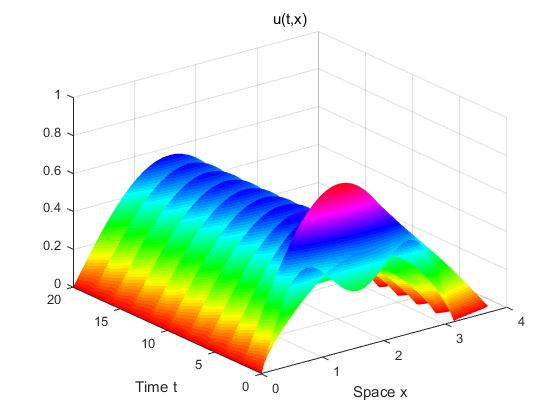}
} }
\subfigure[]{ {
\includegraphics[width=0.30\textwidth]{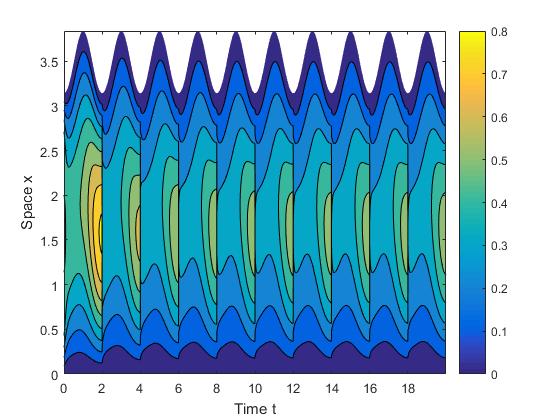}
} }
\subfigure[]{ {
\includegraphics[width=0.30\textwidth]{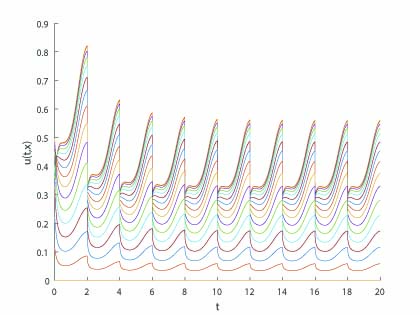}
} }
\caption{\scriptsize $\rho_2(t)=e^{0.1(1-\cos\pi t)}$ and $g(u)$ equals the Ricker function with $b=1.2$ and $r=0.05$.
In this situation, $R_0>1$.  Graphs $(a)-(c)$ indicate that the population approaches a positive periodic steady state. }
\label{tu4}
\end{figure}

Now we consider the case when $g$ is nonmonotone. We fix $g(u)$ as the Ricker function with $b=1.2$ and $r=0.05$ as in (\ref{a03}),
then $g(u)$ satisfies assumptions (A1) and (B3).
\begin{exm}
We fix $g(u)=ue^{0.05-1.2u}$ , then $g'(0)=e^{0.05}$. As Example 4.1, we choose $\rho_1(t)=e^{-0.1(1-\cos\pi t)}$ firstly, then
$$R_0(\rho_1)=\frac{\alpha}{\frac{d\lambda_1}{T}\int_0^T\frac{1}{\rho^2(t)}dt-\frac{1}{T}\ln g'(0)}\approx0.9101<1.$$
One can see from Fig. \ref{tu3} that the population eventually suffers extinction.

We next choose $\rho_2(t)=e^{0.1(1-\cos\pi t)}$. Direct calculations show that
$$R_0(\rho_2)=\frac{\alpha}{\frac{d\lambda_1}{T}\int_0^T\frac{1}{\rho^2(t)}dt-\frac{1}{T}\ln g'(0)}\approx1.1721>1.$$
It is easily seen from Fig. \ref{tu4} that the population approaches a positive periodic steady state.

It can be seen from Example 4.2 that when impulsive harvesting occurs in the form of the Ricker function, the population presents dynamics
similar to the case where the impulsive harvesting occurs in the form of the Beverton-Holt function, that is, the population vanishes in a
periodically evolving habitat with a
small evolution rate and persists in a habitat with a larger evolution rate.
\end{exm}


It is shown in Examples 4.1 and 4.2 that the evolution rate of the domain can impose similar dynamical behaviors of the population,
no matter if $g(u)$ is given as the Beverton-Holt function or the Ricker function. We note that the evolution rate of domain plays an
important role in persistence and extinction of the population, that is, the larger the evolution rate is, the more beneficial it is for
the populations survival. We also note that a large domain evolution rate has a positive effect on the survival of the population when an
impulsive harvesting takes place.

\subsection{The effect of impulsive harvesting}
In order to investigate how impulsive harvesting affects the dynamics of the population in a periodically evolving habitat, we employ numerical
simulations to compare situations when impulsive harvesting does occur or not. We firstly consider the case when monotone impulsive harvesting
takes place, and $g(u)$ is chosen as the Beverton-Holt function.

\begin{figure}[ht]
\centering
\subfigure[]{ {
\includegraphics[width=0.3\textwidth]{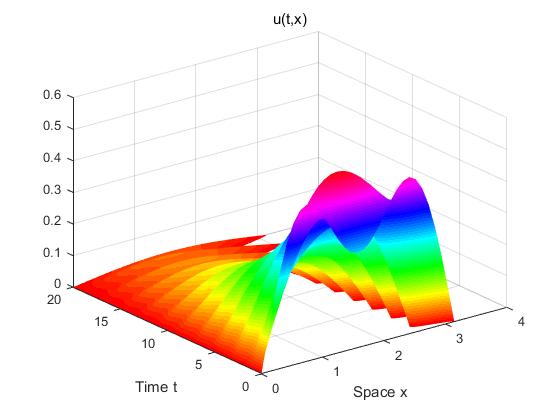}
} }
\subfigure[]{ {
\includegraphics[width=0.3\textwidth]{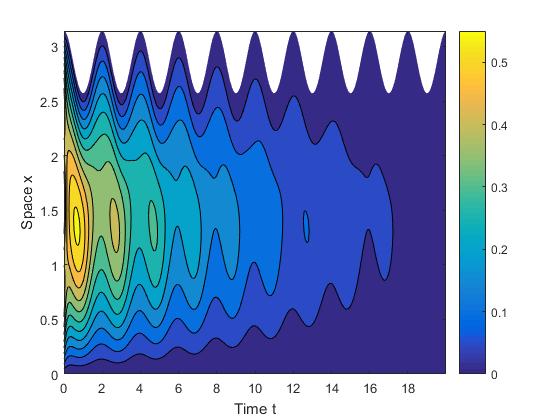}
} }
\subfigure[]{ {
\includegraphics[width=0.3\textwidth]{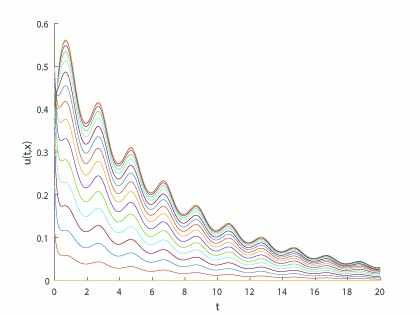}
} }
\caption{\scriptsize  A numerical simulation of the population without impulses and $\rho_1(t)=e^{-0.1(1-\cos\pi t)}$. In this case, $R<1$.
Graphs $(a)-(c)$ show that the population $u(t,x)$ decays to $0$. Graphs $(b)$ and $(c)$ are the cross-sectional view and projection of $u$
on the plane $t-u$, respectively. }
\label{tu5}
\end{figure}

\begin{figure}[ht]
\centering
\subfigure[]{ {
\includegraphics[width=0.30\textwidth]{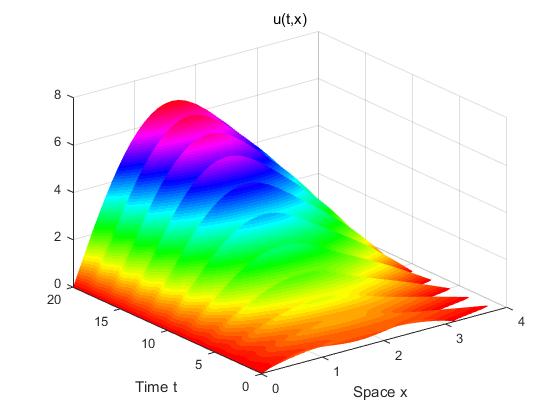}
} }
\subfigure[]{ {
\includegraphics[width=0.30\textwidth]{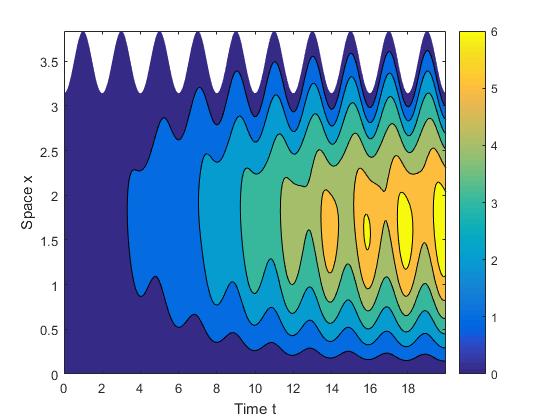}
} }
\subfigure[]{ {
\includegraphics[width=0.30\textwidth]{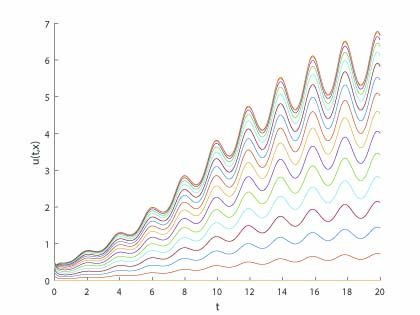}
} }
\caption{\scriptsize  A numerical simulation of the population without impulses and $\rho_2(t)=e^{0.1(1-\cos\pi t)}$. In this situation, $R>1$.
Graphs $(a)-(c)$ show that population $u(t,x)$ evolves to a positive periodic steady state. }
\label{tu6}
\end{figure}

\begin{figure}[ht]
\centering
\subfigure[]{ {
\includegraphics[width=0.30\textwidth]{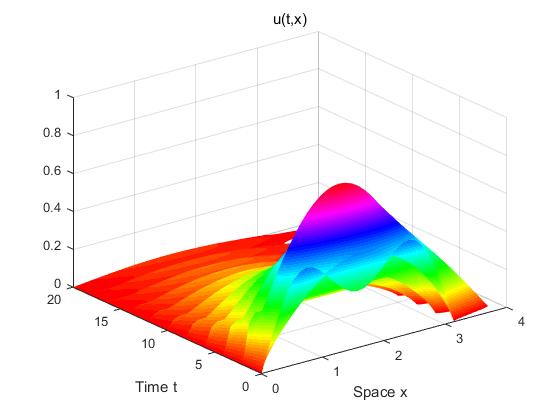}
} }
\subfigure[]{ {
\includegraphics[width=0.30\textwidth]{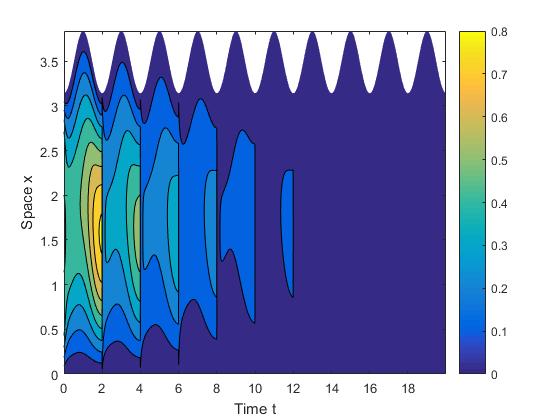}
} }
\subfigure[]{ {
\includegraphics[width=0.30\textwidth]{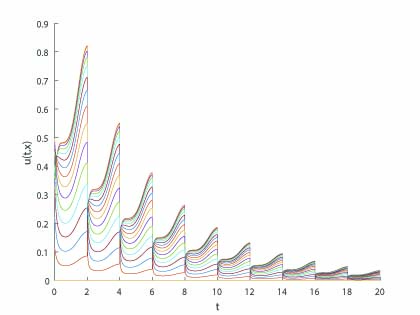}
} }
\caption{\scriptsize In the situation with impulsive harvesting, $\rho_2(t)=e^{0.1(1-\cos\pi t)}$ and $g(u)$ is chosen as the Beverton-Holt
function with $a=10$ and $m=5$, then $R_0<1$. Graph $(a)$ indicates that the population suffers eventual extinction. Graph $(b)$ presents
the case when the domain is periodically evolving. The effect of impulsive harvesting every time $T=2$ can be seen in graph $(c)$, in which the population becomes extinct.  }
\label{tu7}
\end{figure}

\begin{exm}
We first fix $\rho_1(t)=e^{-0.1(1-\cos\pi t)}$. It follows from \eqref{b09} that $R(\rho_1)\approx0.8917<1$ and  impulsive harvesting does not occur.
One can see from Fig. \ref{tu5} that the population suffers extinction eventually. Comparing Fig. \ref{tu1} to Fig. \ref{tu5}, we can conclude
that when impulsive harvesting occurs, population suffers extinction at a faster speed.

Now we fix $\rho_2(t)=e^{0.1(1-\cos\pi t)}$. Now,  without impulsive harvesting, we have $R(\rho_2)\approx1.3302>1.$
It is easily seen from Fig. \ref{tu6} that the population stabilizes to a positive periodic steady state. Then we choose $g(u)=5u/(10+m)$.
One can see from Fig. \ref{tu7} that the population now decays to extinction. We note from Figs. \ref{tu6} and \ref{tu7} that the population
survives in a evolving domain with a large evolution rate, but vanishes when the impulsive harvesting takes place.
\end{exm}


Now we discuss the cases where nonmonotone impulsive harvesting takes place, and $g(u)$ is chosen as the Ricker function.

\begin{figure}[ht]
\centering
\subfigure[]{ {
\includegraphics[width=0.30\textwidth]{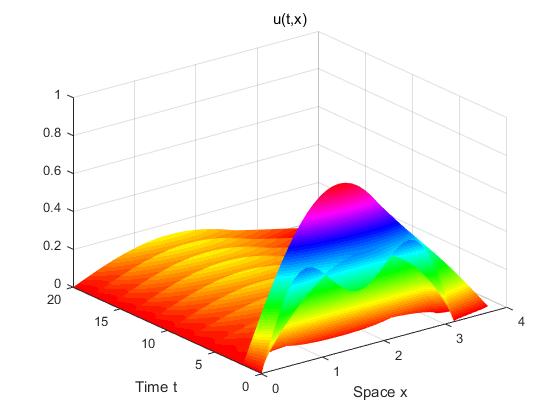}
} }
\subfigure[]{ {
\includegraphics[width=0.30\textwidth]{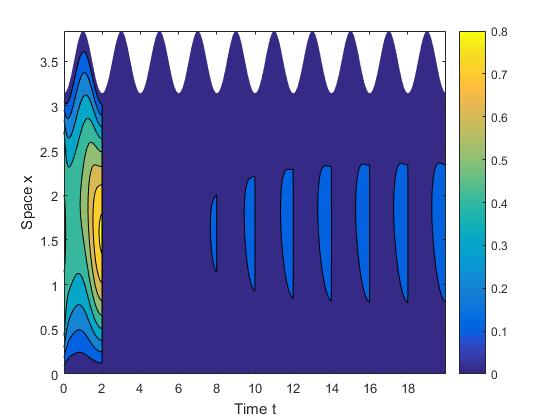}
} }
\subfigure[]{ {
\includegraphics[width=0.30\textwidth]{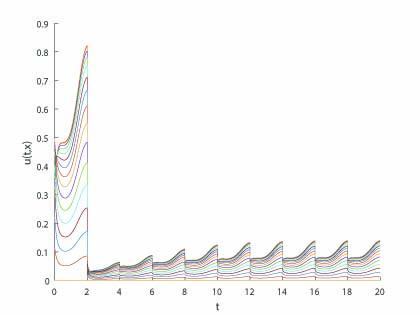}
} }
\caption{\scriptsize In the situation with impulsive harvesting, $\rho_2(t)=e^{0.1(1-\cos\pi t)}$ and $g(u)$ is chosen as the Ricker
function with $b=5$ and $r=0.05$, then $R_0>1$.  Graphs $(a)-(c)$ show that when impulsive harvesting happens, the population drops
sharply at the beginning and evolves to a small positive periodic steady state eventually. }
\label{tu8}
\end{figure}

\begin{exm}
Fix $\rho_1(t)=e^{-0.1(1-\cos\pi t)}$. A comparison of Figs. \ref{tu3} with impulses and \ref{tu5} without impulses reveals that impulsive harvesting
accelerates the extinction of the population.

We next fix $\rho_2(t)=e^{0.1(1-\cos\pi t)}$ and choose $g(u)=ue^{0.05-5u}$. Fig. \ref{tu8} shows that the population approaches a very small positive steady state.
Taken together, Figs. \ref{tu6} without impulses and \ref{tu8} without impulses indicate that the population develops well in an evolving domain
with a large evolution rate, and when impulsive harvesting occurs, the size of the population drops sharply at the beginning but eventually survives.
\end{exm}

Examples 4.3 and 4.4 indicate that when the population lives in a periodically evolving habitat with a small evolution rate,
impulsive harvesting can speed up the extinction of the population no matter which kind of impulse function is chosen. We see that
impulsive harvesting has a negative effect on the survival of the population and can even can lead to the extinction of the population.

\section{\bf Discussion}
A diffusive logistic population model with impulsive harvesting on a periodically evolving domain has been investigated in the present paper.
What we want to study is how the evolution rate of the domain and impulsive harvesting affect on the dynamics of the population.
To address this, we firstly introduce a new threshold value $R_0$ for this impulsive problem; this is the ecological reproduction index,
and is given by an explicit formula and it mainly depends on the evolution rate $\rho(t)$ and the term $g'(0)$, where $g$ is the pulse
function modelling the impulsive harvesting. Then, considering two cases where the function $g(u)$ is either monotone or not,
threshold-type dynamical behaviors of the solution to problem (\ref{a11}) are established. We conclude that when $R_0$ is smaller than 1,
the solution $v(t,y)$ decays to 0 as time goes to infinity no matter which function $g(u)$ is chosen (see Theorems 3.2 and 3.7(i)).
On the other hand, when $R_0$ is greater than 1, it is proved in Theorem 3.6 that the solution $v(t,y)$ converges to a positive
periodic steady state in the monotone case, while in the nonmonotone case, solution $v(t,y)$ converges to an attractive sector (see Theorem 3.7(ii)).

Our numerical simulations further reveal the effects of the evolution rate of the domain and impulsive harvesting on the persistence and
extinction of the population. Examples 4.1 and 4.2 illustrate how the population suffers extinction in a evolving habitat with a small evolution rate,
but survives in one with a large evolution rate.  That is, a large evolution rate can be beneficial for the survival
of the population no matter which kind of impulsive harvesting that occur. Another notable result is that impulsive
harvesting can speed up the population extinction (see Figs. 1, 3 and 5) and has negative effect on the population survival
(see Figs. 6 and 8), and can even lead to the extinction of the population (see Fig. 7).

Our analyses and simulations are based on the one-dimensional case. Recently, Fazly, Lewis and Wang \cite{fm2} have extended the results in \cite{lma} to a higher dimensional habitat without evolving. They provided critical domain parameters and investigated the extinction and persistence of species, depending on the geometry and size of the domain.
We are curious about what new effects the evolution rate will have on the population dynamics when a higher dimensional and evolving region is introduced. Due to the spatial heterogeneity, model parameters, which depend on location, can be considered. Furthermore, a hybrid reaction-advection-diffusion model with a nonlocal discrete time map have been very recently investigated by Fazly, Lewis and Wang in \cite{fm1}. They obtained the existence of travelling wave solutions and provided explicit expressions for the spreading speed. If the nonlocal discrete time map is introduced to our model, it will be interesting to investigate how this will affect the dynamics of the population. This challenges us to a further study.

\end{document}